\documentclass[10pt]{article}
\usepackage{amsfonts}
\usepackage[leqno]{amsmath}
\usepackage{graphicx}
\usepackage{latexsym}
\usepackage{amsmath,amsfonts,amssymb,amsthm,mathrsfs,euscript,makeidx,color}
\usepackage{enumerate}
\usepackage{bbm}
\allowdisplaybreaks[1]
\usepackage{multirow}

\usepackage[numbers, sort&compress]{natbib}

\newcommand{\citedef}[2]{\cite[#2]{#1}}

\def\sqr#1#2{{\vcenter{\vbox{\hrule height.#2pt
				\hbox{\vrule width.#2pt height#1pt \kern#1pt \vrule width.#2pt}
				\hrule height.#2pt}}}}
\def\signed #1{{\unskip\nobreak\hfil\penalty50
		\hskip2em\hbox{}\nobreak\hfil#1
		\parfillskip=0pt \finalhyphendemerits=0 \par}}
\def\endpf{\signed {$\sqr69$}}

\def\3n{\negthinspace \negthinspace \negthinspace }
\def\2n{\negthinspace \negthinspace }
\def\1n{\negthinspace }
\def\bel{\begin{equation}\label}

\def\dbE{\mathbb{E}}
\def\dbF{\mathbb{F}}

\def\dbH{\mathbb{H}}

\def\dbP{\mathbb{P}}

\def\dbR{\mathbb{R}}
\def\dbS{\mathbb{S}}

\def\sU{\mathscr{U}}

\def\sX{\mathscr{X}}

\def\scT{\sc T}
\def\ds{\displaystyle}

\def\ns{\noalign{\ss}}
\def\a{\alpha}

\def\d{\delta}
\def\e{\varepsilon}

\def\l{\lambda}

\def\si{\sigma}
\def\t{\tau}
\def\f{\varphi}

\def\i{\infty}
\def\D{\Delta}
\def\Th{\Theta}

\def\O{\Omega}
\def\cF{{\cal F}}

\def\cK{{\cal K}}

\def\D{\Delta}
\def\Th{\Theta}

\def\O{\Omega}

\def\ss{\smallskip}
\def\ms{\medskip}
\def\bs{\bigskip}
\def\q{\quad}
\def\qq{\qquad}
\def\hb{\hbox}

\def\h1{\outline{$1$}}

\def\hh1{\outline{$1$}}
\def\hh2{\outline{$2$}}
\def\hh3{\outline{$3$}}
\def\hh4{\outline{$4$}}
\def\hh5{\outline{$5$}}
\def\hh6{\outline{$6$}}
\def\hh7{\outline{$7$}}
\def\hh8{\outline{$8$}}
\def\hh9{\outline{$9$}}
\def\hh0{\outline{$0$}}

\def\h{\widehat}
\def\wt{\widetilde}

\def\cd{\cdot}
\def\cds{\cdots}

\def\les{\leqslant}
\def\ges{\geqslant}

\def\({\Big (}
\def\){\Big )}
\def\[{\Big[}
\def\]{\Big]}

\def\bde{\begin{definition}\label}
	\def\ede{\end{definition}}
\def\bel{\begin{equation}\label}
		\def\ee{\end{equation}}
	\def\bt{\begin{theorem}\label}
		\def\et{\end{theorem}}
	\def\bc{\begin{corollary}\label}
		\def\ec{\end{corollary}}
	\def\bl{\begin{lemma}\label}
		\def\el{\end{lemma}}
	\def\bp{\begin{proposition}\label}
		\def\ep{\end{proposition}}
	\def\bex{\begin{example}\label}
		\def\ex{\end{example}}
	\def\bas{\begin{assumption}}
		\def\eas{\end{assumption}}
	\def\br{\begin{remark}\label}
		\def\er{\end{remark}}
	\def\ba{\begin{array}}
		\def\ea{\end{array}}
	\def\ed{\end{document}}

\def\rf{\eqref}

\def\square#1{\vbox{\hrule\hbox{\vrule height#1%
			\kern#1\vrule}\hrule}}
\def\rectangle#1#2{\vbox{\hrule\hbox{\vrule height#1%
			\kern#2\vrule}\hrule}}

\font\tenbb=msbm10 \font\sevenbb=msbm7 \font\fivebb=msbm5

\newfam\bbfam
\scriptscriptfont\bbfam=\fivebb \textfont\bbfam=\tenbb
\scriptfont\bbfam=\sevenbb

\newtheorem{theorem}{Theorem}[section]
\newtheorem{corollary}[theorem]{Corollary}

\newtheorem{lemma}[theorem]{Lemma}
\newtheorem{proposition}[theorem]{Proposition}

\theoremstyle{definition}
\newtheorem{definition}[theorem]{Definition}

\newtheorem{remark}[theorem]{Remark}

\newtheorem{example}[theorem]{Example}

\makeatletter

\@addtoreset{equation}{section}
\makeatother

\usepackage{xcolor}
\makeatletter

\input pdf-trans
\newbox\qbox
\def\usecolor#1{\csname\string\color@#1\endcsname\space}
\newcommand\bordercolor[1]{\colsplit{1}{#1}}
\newcommand\fillcolor[1]{\colsplit{0}{#1}}
\newcommand\outline[1]{\leavevmode%
	\def\maltext{#1}%
	\setbox\qbox=\hbox{\maltext}%
	\boxgs{Q q 2 Tr \thickness\space w \fillcol\space \bordercol\space}{}%
	\copy\qbox%
}
\makeatother
\newcommand\colsplit[2]{\colorlet{tmpcolor}{#2}\edef\tmp{\usecolor{tmpcolor}}%
	\def\tmpB{}\expandafter\colsplithelp\tmp\relax%
	\ifnum0=#1\relax\edef\fillcol{\tmpB}\else\edef\bordercol{\tmpC}\fi}
\def\colsplithelp#1#2 #3\relax{%
	\edef\tmpB{\tmpB#1#2 }%
	\ifnum `#1>`9\relax\def\tmpC{#3}\else\colsplithelp#3\relax\fi
}
\bordercolor{black}
\fillcolor{white}
\def\thickness{.3}

\begin{document}
\title{\bf A Time-Inconsistent
Stochastic Optimal Control Problem in an Infinite Time Horizon}

\author{
Qingmeng Wei\thanks{School of Mathematics and Statistics, Northeast Normal University,
Changchun, Jilin 130024, China  (Email: {\tt weiqm100@nenu.edu.cn}).    This author is supported in part by the National Key R\&D Program of China (No. 2023YFA1009002),   NSF of P. R. China (No. 12371443)  and the Changbai Talent Program of Jilin Province.}~~~~and~~~
Jiongmin Yong\thanks{Department of Mathematics, University of Central Florida, Orlando, FL 32816, USA (Email: {\tt jiongmin.yong@ucf.edu}).                          This author is supported in part by NSF Grant DMS-2305475.}
                         ~~~~
                         }

\maketitle

\ms

\centerline{\it This paper is dedicated to Professor Shige Peng on the occasion of his 80th birthday}

\bs
\begin{abstract} This paper is concerned with a time-inconsistent
stochastic optimal control problem in an infinite
time horizon with a non-degenerate diffusion in the state equation. A major assumption is that people become rational after a large time. Under such a condition, the problem in an infinite time horizon can be decomposed into two parts: a non-autonomous time-consistent problem in an infinite time horizon and a time-inconsistent problem in a finite time horizon. Then an equilibrium strategy will be constructed. Both Bolza type problem and recursive cost problem are considered.

\end{abstract}


\ms


\section{Introduction --- Time-Inconsistency.}

Let $(\O,\cF,\dbF,\dbP)$ be a complete filtered probability space satisfying the usual condition, on which a $d$-dimensional standard Brownian motion $W(\cd)\equiv(W_1(\cd),\cds,W_d(\cd))$ is defined, %\footnote{If the dimension of the Brownian motion is $d>1$, our approach still works, with a little more complicated notations.}
whose natural filtration is $\dbF=\{\cF_t\}_{t\ges0}$
(augmented by all the $\dbP$-null sets). Denote $\dbR_+=[0,\i)$.
%
%$$L^2_{\cF_t}(\O;\dbR^n)=\Big\{\xi:\O\to\dbR^n\bigm|\xi\hb{ is $\cF_t$-measurable, }\dbE|\xi|^2<\i\Big\}.$$
%
We consider the following controlled stochastic differential equation (SDE, for short) with $0\les t<T\les\i$:
\bel{state1}\left\{\2n\ba{ll}
\ds dX(s)=b(s,X(s),u(s))ds+\sum_{i=1}^d\si_i(s,X(s),u(s))dW_i(s),\qq s\in[t,T)],\\
\ns\ds X(t)=x,\ea\right.\ee
where $b,\si_1,\cds,\si_d:\dbR_+\times\dbR^n\times U\to\dbR^n$ are continuous ($U\subseteq\dbR^m$ is a closed set), and $(t,x)\in[0,T)]\times\dbR^n$ is called the {\it initial pair}. The set $[t,T)]$ stands for $[t,T]$ if $T<\i$ and $[t,\i)$ if $T=\i$. The pair $(X(\cd),u(\cd))$, with $X:[t,\i)\times\O\to\dbR^n$, called the {\it state process}, and $u:[t,\i)\times\O\to U$, called the {\it control process}, is called a {\it state-control pair} on $[t,T)]$. Now, let us introduce the set
\bel{cU}\sU[t,T]=\Big\{u:[t,T)]\times\O\to U\bigm|u(\cd)\hb{ is
$\dbF$-progressively measurable}\Big\},\qq0\les t<T\les\i.\ee
Any $u(\cd)\in\sU[t,T]$ is called an {\it admissible control} on  $[t,T)]$. For any initial pair $(t,x)\in[0,T)]\times\dbR^n$ and $u(\cd)\in\sU[t,\i]$, under some mild conditions on $b(\cd)$ and $\si(\cd)\equiv(\si_1(\cd),\cds,\si_d(\cd))$, state equation (\ref{state1}) admits a unique strong solution  $X(\cd)\equiv X(\cd\,;t,x,u(\cd))$, which is in the following set
\bel{X}\ba{ll}
\ns\ds\sX[t,\i]=\Big\{X:[t,\i)\times\O\to\dbR^n\bigm|X(\cd)\hb{ is $\dbF$-progressively measurable,}\\
\ns\ds\qq\qq\qq\qq\qq\qq\qq\qq\dbE\[\sup_{s\in[t,T]}|X(s)|^2\]<\i,~\forall t<T<\i\Big\} .\ea\ee
To measure the performance of the control process $u(\cd)\in\sU[t,T]$, we introduce the following cost functional on $[t,T)]$:
\bel{J^T}J^T_*(t,x;u(\cd))=\dbE\[\int_t^T\l(s-t)g(s,X(s),u(s))ds
+\l(T-t)h(X(T))\],\ee
where $g:\dbR_+\times\dbR^n\times\dbR^m\to\dbR_+$, $h:\dbR^n\to\dbR_+$ are continuous and bounded, $\l:\dbR_+\to\dbR_+$ is continuous and strictly decreasing with
\bel{t(s)}\l(0)=1,\q\lim_{s\to\i}\l(s)=\l(\i)=0,\q\l(\cd)\in L^1(0,\i;\dbR).\ee
It suffices to assume that $g(\cd\,,\cd\,,\cd)$ and $h(\cd)$ are bounded from below. Letting them be non-negative is just for convenience. We have seen that the dependence on $h(\cd)$ is emphasized by a subscript $*$ (since the expression of $h(\cd)$ might be complicated later), which will be useful later. In the case $T=\i$, the term involving $h(\cd)$ on the right-hand side of \rf{J^T} automatically vanishes (since $h(\cd)$ is assumed to be bounded), and $J^\i(t,x;u(\cd))$ will be written instead of $J^\i_*(t,x;u(\cd))$. The role played by $\l(\cd)$ is to discount the running cost rate function on $[t,T)]$. Thus, it is meaningless to consider the cost functional for the past. Hence, if the current time is $t$, the following cost functional is meaningless:
$$\dbE\[\int_{\wt t}^T\l(s-t)g(s,X(s),u(s))ds\],\qq0\les\wt t<t<T\les\i.$$
Or, we might trivially extend $\l(r)$ to be zero for $r<0$. Now we pose the following optimal control problem.

\ms

\bf Problem (P)$^T_*$. \rm For given $t\in[0,T)]$ and $x\in\dbR^n$, find a $\bar u^{t,T}(\cd\,;x)\in\sU[t,T]$ such that
\bel{V=inf}J^T_*(t,x;\bar
u^{t,T}(\cd\,;x)))=\inf_{u(\cd)\in\sU[t,T]}J^T_*(t,x;u(\cd))\equiv \bar V^{t,T}_*(x).\ee

\ss

Any $\bar u^{t,T}(\cd\,;x)\in\sU[t,T]$ satisfying \rf{V=inf} is called a {\it pre-committed open-loop optimal control} of Problem (P)$^T_*$ for the initial state $x\in\dbR^n$, at $t\in[0,T)]$, the
corresponding state process $\bar X^{t,T}(\cd\,;x)\equiv X(\cd\,;t,x,\bar u^{t,T}(\cd\,;x))$ and the pair $(\bar X^{t,T}(\cd\,;x),\bar u^{t,T}(\cd\,;x))$ are called the corresponding {\it pre-committed open-loop optimal state process} and {\it pre-committed open-loop optimal pair} of Problem (P)$^T_*$, at $t\in[0,T)]$, respectively. The function $\bar V^{t,T}_*(\cd)$ defined by \rf{V=inf} is called the {\it pre-committed value function} of Problem (P)$^T_*$ at $t\in[0,T)]$. Here, the word ``pre-committed'' is used, which means that based on the initial pair $(t,x)\in[0,T)]\times\dbR^n$, once the optimal control $\bar u^{t,T}(\cd)$ is found based on this initial pair and applied on $[t,T)]$, regardless what happens later, the control (which might not be optimal anymore) will not be changed (pre-committed).

\ms

On the other hand, since $\bar u^{t,x}(\cd)$ might not be optimal when $(t,x)$ becomes $(t',X(t';t,x,\bar u^{t,x}(\cd)))$ (for $t'>t$), one might want to pose the following problem.

\ms

\bf Problem (N)$^T_*$. \rm For any initial pair $(t,x)\in[0,T)]\times\dbR^n$, find a $\bar u^T(\cd\,;t,x)\in\sU[t,T]$ such that
\bel{V2}J^T_*(t,x;\bar
u^T(\cd\,;t,x))=\inf_{u(\cd)\in\sU[t,T]}J^T_*(t,x;u(\cd))\equiv \bar V^T_*(t,x).\ee
Namely, for any initial pair $(t,x)\in[0,T)]\times\dbR^n$, define
$$(\bar X^T(s;t,x),\bar u^T(s;t,x))=(\bar X^{t,T}(s;x),\bar u^{t,T}(s;x)),\qq s\in[t,T)].$$

Any $\bar u^T(\cd\,;t,x)\in\sU[t,T]$ satisfying \rf{V2} is called a {\it naive open-loop optimal control} of Problem (N)$^T_*$ for the initial pair $(t,x)\in[0,T)]\times\dbR^n$, the
corresponding state process $\bar X^T(\cd\,;t,x)\equiv X(\cd\,;t,x,\bar u^T(\cd\,;t,x))$ and the pair $(\bar X^T(\cd\,;t,x),\bar u^T(\cd\,;t,x))$ are called the corresponding {\it naive open-loop optimal state process} and {\it naive open-loop optimal pair}, of Problem (N)$^T_*$, respectively. The function $\bar V^T_*(\cd\,,\cd)$ defined by \rf{V2} is called the {\it naive value function} of Problem (N)$^T_*$. In what follows, for simplicity, when $T=\i$, we call Problems (P)$^\i$ and (N)$^\i$, respectively instead, and denote the value functions by
$\bar V^{\sc t,\i}(x)$ and $\bar V^\i(t,x)$, respectively.

\ms

In the above, the term ``naive'' was used for the following reason. Control $\bar u^T(\cd\,;t,x)$ is optimal only at $(t,x)\in[0,T)]\times\dbR^n$. When the time becomes $t'>t$, in order the performance remains optimal at $t'$, the control has to be taken as $\bar u^T(\cd\,;t',\bar X(t'))$ (selected based on the pair $(t',\bar X(t'))$), since $\bar u^T(\cd\,;t,x)$ might not be optimal on $[t',T)]$ anymore. Thus, as time goes by, the controller has to change the value of the control instantaneously, in order to keep the performance to be constantly optimal, which is practically impossible.

\ms

The simple difference between the pre-committed optimal control $\bar u^{t,T}(\cd\,;x)$ and naive optimal control $\bar u^T(\cd\,;t,x)$ is that the former is a function of $x$ with $t$ being a parameter. Whereas, the latter is a function of $(t,x)$. From the actions, we clearly see the difference between these two.

\ms

Mathematically, the time-consistency is described by the fact whether the optimal control contains just one running time or two. To get more feeling about this, let us make an observation. Note
\bel{wt J}\ba{ll}
\ns\ds J^\i(t,x;u(\cd))\equiv\dbE\[\int_t^\i\l(s-t)g(s,X(s),u(s))ds\]\\
\ns\ds\qq\qq\qq=\dbE\[\int_0^\i\l(s)g(s+t,X(s+t),u(s+t))ds\]
\equiv\dbE\[\int_0^\i\l(s)\wt g(s,\wt X(s),\wt u(s))ds\],\ea\ee
where
$$\wt g(s,x,u)=g(s+t,x,u),$$
and
$(\wt X(\cd),\wt u(\cd))$ is the state-control pair of the following controlled SDE:
\bel{wt X}\left\{\2n\ba{ll}
\ds d\wt X(s)=\wt b(s,\wt X(s),\wt u(s))ds+\sum_{i=1}^d\wt\si_i(s,\wt X(s),\wt u(s))d\wt W^t_i(s),\qq s\ges0,\\
\ns\ds\wt X(0)=x,\,\ea\right.\ee
with
$$\wt b(s,x,u)=b(s+t,x,u),\q\wt\si_i(s,x,u)=\si_i(s+t,x,u),\q s\ges0,$$
and $\wt W^t(s)=W(s+t)$ being the Brownian motion starting from a  random vector $W(t)$ (instead of 0) with a normal distribution. Then, $(X(\cd),u(\cd))$ and $(\wt X(\cd),\wt u(\cd))$ have the same distribution. But as processes, it holds
\bel{ne}(\wt X(\cd),\wt u(\cd))\ne(X(\cd),u(\cd)),\qq(\wt X(s-t),\wt u(s-t))=(X(s),u(s)),\q s\in[t,\i).\ee
The above suggests us considering a classical non-autonomous stochastic optimal control problem in an infinite horizon. More precisely,  we introduce the following state equation:
\bel{hat X}\left\{\2n\ba{ll}
\ds d\h X(s)=\h b(s,\h X(s),\h u(s))ds+\sum_{i=1}^d\h\si_i(s,\h X(s),\h u(s))d\h W_i(s),\qq s\in[t,\i),\\
\ns\ds\h X(t)=x,\,\ea\right.\ee
with $\h b(\cd),\h\si(\cd)$ being given and
$\h W(\cd)$ being a $d$-dimensional Brownian motion starting from a random vector having a normal distribution, and with cost functional (with $\h X(\cd)$ being the corresponding state process)
\bel{hat J}\h J^\i(t,x;\h u(\cd))=\dbE\int_t^\i\h g(s,\h X(s),\h u(s))ds,\ee
where $\h g:\dbR_+\times\dbR^n\times\dbR^m\to\dbR_+$ is decay so that the cost functional is well-defined. We formulate the following time-consistent optimal control problem.
\ms

{\bf Problem ($\h{\rm C}$)$^\i$.} For $(t,x)\in\dbR_+\times\dbR^n$, find $\h u^\i(\cd)\equiv\h u^\i(\cd\,;t,x)\in\sU[t,\i]$ such that \rf{hat X} is satisfied and
$$\h J^\i(t,x;\h u^\i(\cd))=\inf_{\h u(\cd)\in\sU[t,\i]}\h J^\i(t,x;\h u(\cd))=\h V^\i(t,x).$$

This problem will be solved in the next section. Note that $\h V^\i(t,x)$ is not explicitly depending on the Brownian motion, because it can be characterized as the unique viscosity solution of the corresponding HJB equation (see \cite{Yong-Zhou 1999}). Hence, by taking
$$\h b(s,x,u)=\wt b(s,x,u),\ \h \si(s,x,u)=\wt\si(s,x,u),\ \h g(s,x,u)=\l (s)\wt g(s,x,u),\ \h W(s)=\wt W^t(s),$$
we see that \rf{hat X} coincides with \rf{wt X} and \rf{hat J} coincides with \rf{wt J}. Therefore, from \rf{ne}, for given $(t,x)\in\dbR_+\times\dbR^n$, naive open-loop optimal pair satisfies
\bel{X,u}(\bar X^\i(s;t,x),\bar u^\i(s;t,x))=(\h X^\i(s-t;t,x),\h u^\i(s-t;t,x)),\qq s\in[t,\i).\ee
the value $\bar V^\i(s,\bar X^\i(s;t,x))$ of the value function $\bar V^\i(\cd\,,\cd)$ of Problem (N)$^\i$ at $s>t$ satisfies
\bel{}\bar V^\i(s,\bar X^\i(s;t,x))=\h V^\i(s-t,x),\qq s\in[t,\i).\ee
From the above, we see how the two running times ($s$ and $t$) in $(\bar X^\i(s;t,x),\bar u^\i(s;t,x))$ and $V^\i(s,\bar X^\i(s;t,x))$ might appear, and how these quantities depends on the initial time $t$, which exhibits the so-called {\it time-inconsistency}.

\ms

A natural question is this: Is it possible to have a certain discount function so that a naive (or a pre-committed) optimal control is {\it time-consistent}, or intuitively, is the optimal control only depending on single running time $s$ (so that one does not have to change the control instantaneously as time goes by)? This is possible if
\bel{l0}\l(\t)=e^{-\d\t},\qq\t\ges0,\ee
for some $\d>0$, called {\it discount rate}. We call such a function $\l(\cd)$ an {\it exponential discount function}. In this case, for $0<T\les\i$, the cost functional can be written as (compare \rf{J^T})
\bel{J^i2}J^{\d,T}_*(t,x;u(\cd))=\dbE\[\int_t^Te^{-\d(t-s)}g(s,X(s),u(s))ds
+e^{-\d(T-t)}h(X(T))\],\ee
indicating the dependence on $\d>0$. Then, minimizing $J^{\d,T}_*(t,x;u(\cd))$ is equivalent to minimizing
\bel{cost2}e^{-\d t}J^{\d,T}_*(t,x;u(\cd))=\dbE\[\int_t^T
e^{-\d s}g(s,X(s),u(s))ds+e^{-\d T}h(X(T))\]
\equiv\h J_*^{\,\d,T}(t,x;u(\cd)),\ee
which, with bounded $g(\cd\,,\cd\,,\cd)$ and $h(\cd)$, is a classical Bolza type cost functional (in $[t,T)]$). Hence,
\bel{V_0}\h V^{\d,T}_*(t,x)\equiv\inf_{u(\cd)\in\sU[t,T]}
\h J^{\,\d,T}_*(t,x;u(\cd))=e^{-\d t}V^{\d,T}_*(t,x).\ee
It is known that the dynamic programming principle (DPP, for short) holds for the current classical Bolza problem (because the proof of DPP is a local argument). Thus, for all $t'>t\ges0$,
\bel{DPP1}\h V_*^{\d,T}(t,x)=\inf_{u(\cd)\in\sU[t,t']}\dbE\[\int_t^{t'}e^{-\d s}g(s,X(s),u(s))ds+\h V_*^{\d,T}(t',X(t'))\].\ee
By directly considering $J^{\d,T}_*(t,x;u(\cd))$ instead of $e^{-\d t}\h J^{\,\d,T}_*(t,x;u(\cd))$, we have the corresponding DPP:
$$\ba{ll}
\ns\ds V^{\d,T}_*(t,x)=e^{\d t}\h V^{\d,T}_*(t,x)
=e^{\d t}\inf_{u(\cd)\in\sU[t,t']}\dbE\[\int_t^{t'}e^{-\d s}
g(s,X(s),u(s))ds+\h V^{\d,T}_*(t',X(t'))\]\\
\ns\ds=\inf_{u(\cd)\in\sU[t,t']}\dbE\[\int_t^{t'}e^{-\d(s-t)}
g(s,X(s),u(s))ds+e^{\d t}\h V^{\d,T}_*(t',X(t'))\]\\
\ns\ds=\inf_{u(\cd)\in\sU[t,t']}\dbE\[\int_t^{t'}e^{-\d(s-t)}
g(s,X(s),u(s))ds+e^{-\d(t'-t)}V^{\d,T}_*(t',X(t'))\].\ea$$
Hence, if $\bar u^T(\cd\,;t,x)\in\sU[t,T]$ is (a naive) an open-loop optimal control for the initial pair $(t,x)\in[0,T)]\times\dbR^n$, then it holds (for $t<t'<T$)
$$\ba{ll}
\ns\ds V^{\d,T}_*(t,x)=J^{\d,T}_*(t,x;\bar u^T(\cd\,;t,x))\\
\ns\ds=\dbE\[\int_t^{t'}e^{-\d(s-t)}g(s,\bar X^T(s),\bar u^T(s))ds
+e^{-\d(t'-t)}J^{\d,T}_*(t',\bar X^T(t');\bar u^T(\cd)\big|_{[t',T)]})\]\\
\ns\ds\ges\dbE\[\int_t^{t'}e^{-\d(s-t)}g(s,\bar X(s),\bar u(s))ds+e^{-\d(t'-t)}\bar V^{\d,T}_*(t',\bar X(t'))\]\\
\ns\ds\ges\inf_{u(\cd)\in\sU[t,t']}\dbE\[\int_t^{t'}e^{-\d(s-t)}
g(s,X(s),u(s))ds+e^{-\d(t'-t)}V_*^{\d,T}(t',X(t'))\]=V^{\d,T}_*(t,x).\ea$$
Consequently, all the equalities in the above hold, in particular,
$$J^{\d,T}_*(t',\bar X^T(t');\bar u^T(\cd)\big|_{[t',T)]})=V^{\d,T}_*(t',\bar X^T(t')).$$
This means that $\bar u^T(\cd\,;t,x)\big|_{[t',T)]}$, an open-loop optimal control selected for the initial pair $(t,x)$, remains an open-loop optimal on $[t',T)]$ for any $T>t'>t\ges0$. This is exactly the meaning of {\it time-consistency}. Clearly, people desire such a property for optimal controls. However, in reality, the time-consistency (the above-mentioned feature) could be lost, due to people's time preferences. Normal people usually over weight the immediate utility, which can be mathematically described by the {\it general discounting}, including the so-called {\it non-exponential discounting}, or {\it hyperbolic discounting}. It turns out that if we do not assume \rf{l0}, then Problem (N)$^T_*$ (with $0<T\les\i$) is {\it time-inconsistent}, in general. Namely, a naive open-loop optimal control found for the initial pair $(t,x)\in[0,T)]\times\dbR^n$ will not stay optimal thereafter (it will be staying naive later on). The purpose of solving the problem is to find an equilibrium strategy/control, if it is possible, which has the following characteristics: If $t\in[0,T)]$ is the current time and such a (locally optimal) control $\bar u^T(\cd\,;t,x)$ is applied, then on the rest $(t,T)]$ of this interval, this control will not be changed and it will stay (locally) optimal. Mathematically, the found control is to be proved as a function of single running time variable $s\in[t,T)]$, showing its time-consistency (and local optimality). We emphasize that the equilibrium strategy/control one is looking for is depending on the initial pair $(t,x)\in[0,T)]\times\dbR^n$, is defined on $[t,T)]$ and is staying (locally) optimal in the rest of the interval $(t,T)]$.

\ms

For finite time horizon problems, the study of time-inconsistency with time-consistent equilibrium solutions have been carried out by many authors, for the case that the diffusion of the state equation being non-degenerate, we mention \cite{Yong 2012}, in which the equilibrium HJB equation was derived the first time, see also \cite{Wei-Yong-Yu 2017, Wei-Yu 2018, Wang-Yong-Zhou 2024} and the references cited therein for the follow-up works. For the case of linear-quadratic problems, see \cite{Hu-Jin-Zhou 2012, Yong 2017, Hu-Jin-Zhou 2017, Lu-Ma 2024, Alia-Alia 2024} and the references cited therein.

\ms

For the problems in an infinite time horizon, if one tries to use the usual differential game approach backwardly (see \cite{Pollak 1968,Yong 2012}), then one faces the difficulty that final player does not exist. Fortunately, there is a  remedy for this. We adopt the following natural assumption: After a long time, normal people will gradually become rational. This is true from many points of view. Thus, we accept that
\bel{l1}\l(\t)=e^{-\d\t},\qq \t\ges T_0,\ee
for some large $T_0>0$. Such a condition can be used to overcome the difficulty mentioned earlier. In this paper, we are going to use the condition \rf{l1}.

\ms

In what follows, we assume that all the later selves are going to do their best. Our goal is to investigate infinite time horizon Problem (N)$^\i$ under condition \rf{l1}. Note that in our case, for any given $(t,x)\in\dbR_+\times\dbR^n$, we have
\bel{J=J+V}\ba{ll}
\ns\ds J^\i(t,x;u(\cd))=\dbE\[\int_t^{t+T_0}\l(s-t)g(s,X(s),u(s))ds+
\int_{t+T_0}^\i e^{-\d(s-t)}g(s,X(s),u(s))ds\]\\
\ns\ds\qq\qq\qq=\dbE\[\int_t^{t+T_0}\l(s-t)g(s,X(s),u(s))ds+e^{\d t}
\int_{t+T_0}^\i e^{-\d s}g(s,X(s),u(s))ds\]\\
\ns\ds\qq\qq\qq=\dbE\1n\int_t^{t+T_0}\3n\l(s-t)g(s,X(s),u(s))ds\1n+\1n e^{\d t}\h J^{\,\d,\i}(t\1n+\1n T_0,X(t\1n+\1n T_0);u(\cd+\1n t\1n+\1n T_0)).\ea\ee
Then the controller will solve a classical non-autonomous stochastic optimal control on $[t+T_0,\i)$ to have
\bel{hV}\inf_{u(\cd)\in\sU[t+T_0,\i]}\h J^{\,\d,\i}(t+T_0,X(t+T_0);u(\cd))=\h V^{\d,\sc\i}(t+T_0,X(t+T_0)).\ee
Next, on $[t,t+T_0]$, the controller solves a time-inconsistent stochastic optimal control problem (in a finite time horizon) with the cost functional being:
\bel{J=V+hJ}\dbE\[\int_t^{t+T_0}\l(s-t)g(s,X(s),u(s))ds+e^{\d t }\h V^{\d,\i}(t+T_0,X(t+T_0))\].\ee
Then, the so-called equilibrium value function is given by
\bel{V}\bar V^\i(t,x)=\dbE\[\int_t^{t+T_0}\l(s-t)g(s,\bar X(s),\bar u(s))ds+e^{\d t}\h V^{\d,\i}(t+T_0,\bar X(t+T_0))\],\ee
where $(\bar X(\cd),\bar u(\cd))$ is an equilibrium pair of Problem (N)$^{t,t+T_0}_*$, the time-inconsistent problem on a finite time interval $[t,t+T_0]$.

\ms

This provides the basic idea of solving our time-inconsistent stochastic optimal control problem in an infinite time horizon. We will construct the corresponding time-consistent equilibrium strategy as well as the characteristics for the equilibrium value function.

\ms

Time-inconsistent problem is very old. One can find early literatures from 18 century (see \cite{Hume 1739,Smith 1759}). The first work using mathematics was due to Strotz \cite{Strotz 1955}, followed by some other authors \cite{Pollak 1968, Miller-Salmon 1985, Tesfatsion 1986, Laibson 1997, Karp-Lee 2003, Herings-Rohde 2006, Caplin-Leahy 2006, Grenadier-Wang 2007, Ekeland-Lazrak 2010}, to mention a few. The survey paper by Palacios-Huetta \cite{Palacios-Huerta 2003} is very useful for those who are interested in the history of time-inconsistent problems.

\ms

The rest of the paper is organized as follows. In Section 2, we present relevant results of non-autonomous optimal control problem. Our approach is a little different from that found in \cite{Baumeister-Leitao-Silva 2007}. Section 3 gives the  construction of an equilibrium strategy over an infinite time horizon. In Section 4, we look at the problem in an infinite time horizon with recursive cost functional. Finally, some conclusion comments are collected in Section 5.

\section{A Non-Autonomous Classical Problem in Infinite Horizon}\label{Sec:non-autonomous}

In what follows, we make some conventions. Any vector is assumed to be a column vector. For differentiable function $f:\dbR^n\to\dbR$, its gradient $f_x(x)$ is a row vector, its Hessian $f_{xx}(x)$ is an $(n\times n)$ symmetric matrices which is denoted by $\dbS^n$, and
\bel{D^*}\D^*[0,T]=\{(\t,t)\in[0,T]\times[0,T]\bigm|0\les \t\les t\les T\},\qq\forall0<T\les\i.\ee
Finally, $K>0$ will be a generic constant which could be different from line to line.

\ms

In this section, we are looking at a non-autonomous stochastic optimal control problem. See \cite{Baumeister-Leitao-Silva 2007}, for a  different approach. The state equation is \rf{hat X} with the cost functional \rf{hat J}. We introduce the following assumption for the coefficients of the state equation.

\ms

{\bf($\h{\bf H}$1)} The maps $\h b,\h\si_1,\cds,\h\si_d:\dbR_+\times\dbR^n\times U\to\dbR^n$ are continuous, regular enough, and for some constant $L>0$,
\bel{|b|}|\h b(t,0,u)|+\sum_{i=1}^d|\h\si_i(t,0,u)|\les L,\qq(t,u)\in\dbR_+\times U,\ee
and
\bel{|b-b|}|\h b(t,x,u)-\h b(t,x',u)|+\sum_{i=1}^d|\h\si_i(t,x,u)-\h\si_i(t,x',u)|\les L|x-x'|, \q(t,u)\in\dbR_+\times U,~x,x'\in\dbR^n.\ee
Moreover, the following non-degeneracy condition holds: For some $\e>0$,
\bel{sisi>d}\h\si(t,x,u)\h\si(t,x,u)^\top\ges\e I_n,\q\forall(t,x,u)\in\dbR_+\times\dbR^n\times U.\ee

\ms

The above \rf{sisi>d} implies
\bel{d>n}d\ges n=\hb{rank}\,\si(t,x,u),\qq(t,x,u)\in\dbR_+\times\dbR^n\times U.\ee
For the running cost rate function, we introduce the following.

\ms

{\bf($\h{\bf H}$2)} The map $\h g:\dbR_+\times\dbR^n\times U\to\dbR_+$ is continuous, regular enough, and for some $\h\f(\cd)\in L^1(\dbR_+;\dbR)$, it holds that
\bel{h-g}0\les\h g(t,x,u)\les\h\f(t),\qq(t,x,u)\in\dbR_+\times\dbR^n\times U.\ee

\ms

Under ($\h{\rm H}$1)--($\h{\rm H}$2), the infinite horizon optimal control problem, i.e., Problem $(\h{\rm C})^\i$ presented in the previous section, is well-formulated. The regularity of the coefficients assumed in $(\h{\rm H}1$)--($\h{\rm H}2$) is just for convenience, and can be relaxed. To solve such a problem, we introduce its finite horizon version. That is to consider the following cost functional ($t<T<\i$)
$$\h J^{\,\scT}(t,x;\h u(\cd))=\dbE\int_t^T\h g(s,\h X(s),\h u(s))ds,$$
and the following optimal control problem in a finite time horizon $[0,T]$:

\ms

{\bf Problem ($\h{\rm C}$)$^T$.} For $(t,x)\in[0,T]\times\dbR^n$, find $\h u^T(\cd)\in\sU[t,T]$ such that
$$\h J^{\,\scT}(t,x;\h u^T(\cd))=\inf_{\h u(\cd)\in\sU[t,T]}\h J^{\,\scT}(t,x;\h u(\cd))=\h V^{\scT}(t,x).$$
In the above,  $\h u^T(\cd)\in\sU[t,T]$ is called an open-loop optimal control of Problem ($\h{\rm C}$)$^T$ for the initial pair $(t,x)\in[0,T]\times\dbR^n$; $\h V^{\scT}(\cd,\cd)$ is called the value function of  Problem ($\h{\rm C}$)$^T$.

\ms

For $(t,x,u,p,P)\in[0,T]\times\dbR^n\times U\times\dbR^n\times\dbS^n$, we set
$$
 \h\dbH(t,x,u,p^\top\2n,P)=p\h
b(t,x,u)+{1\over2}\sum_{i=1}^d\h\si_i(t,x,u)^\top
P\h\si_i(t,x,u)+\h g(t,x,u).  $$
Moreover, suppose (similar to \cite{Yong 2012}) there exists a map $\h\psi : \dbR_+\times\dbR^n \times\dbR^n\times\dbS^n \to U$ with needed
regularity such that
\bel{hpsi}\h\psi(t,x,p^\top\2n,P)\in\arg\min_{u\in U}\h\dbH(t,x,u,p^\top\2n,P),\q (  t, x, p, P) \in\dbR_+\times\dbR^n \times\dbR^n\times\dbS^n .\ee

For Problem ($\h{\rm C})^T$, we have the following result.

\bl{monotone} \sl Let {\rm($\h{\rm H}$1)--($\h{\rm H}$2)} and \rf{hpsi} hold.

\ms

{\rm(i)} For given $T>0$, the value function $\h V^T(\cd\,,\cd)$ of Problem {\rm($\h{\rm C})^T$} is the unique classical solution to the following HJB equation
\bel{HJB4}\left\{\2n\ba{ll}
\ds\h V^T_t(t,x)+\h\dbH\big (t,x, \h\psi(t,x,V^T_x(t,x),V^T_{xx}(t,x)),\h V^T_x(t,x),\h V^T_{xx}(t,x)\big)=0,\q(t,x)\in[0,T]\times\dbR^n,\\
\ns\ds\h V^T(T,x)=0,\qq x\in\dbR^n.\ea\right.\ee
Moreover, Problem {\rm($\h{\rm C})^T$} admits the unique optimal control $\h u^T(\cd)\in\sU[t,T]$, which can be represented by
\bel{}\h u^T(s)=\h\psi\big(s,\h X^{t,T}(s),\h V^{T}_x(s,\h X^{t,T}(s)),\h V^{T}_{xx}(s,\h X^{t,T}(s))\big),\qq s\in[t,T].\ee

{\rm(ii)} Further, there exists a function $\h V^\i(\cd\,,\cd)\in C^{1,2}([0,T]\times\dbR^n)$ such that
\bel{V to V}\lim_{T\to\i}\h V^T(t,x)=\h V^\i(t,x)\ee
uniformly in $(t,x)\in\dbR_+\times\dbR^n$. Moreover, $\h V^\i(\cd\,,\cd)$ is the unique classical solution of the following HJB equation
\bel{HJB1}\left\{\2n\ba{ll}
\ds\h V_t^{\sc\i}(t,x)+\h\dbH\big(t,x,\h\psi(t,x,\h V^\i_x(t,x),\h V_{xx}^\i(t,x)),\h V^{\sc\i}_x(t,x),\h V^{\sc\i}_{xx}(t,x)\big)=0,\qq(t,x)\in\dbR_+\times\dbR^n,\\
\ns\ds\h V^\i(\i,x)\equiv\lim_{T\to\i}\h V^{\sc\i}(T,x)=0,\qq x\in\dbR^n.\ea\right.\ee

\el

\it Proof. \rm  The assertion  (i) is classical in control theory. The readers are referred to \cite{Yong-Zhou 1999}.
We now focus on (ii).

\ms

{\it Step 1: Convergence of $\h V^T(\cd\,,\cd)$.}
Since $\h g(\cd)$ satisfies \rf{h-g}, for all $(t,x)\in[0,T]\times\dbR^n$ and $\h u(\cd)\in \sU[t,T]$, we see that
$\ds 0\les\h J^{\,\scT}(t,x;\h  u(\cd))\les K,$
which leads to
$$0\les\h V^T(t,x)\les K,\qq(t,x)\in[0,T]\times\dbR^n.$$
Now, for $T<T'$, consider
\bel{HJB'}\left\{\2n\ba{ll}
\ds\h V^{T'}_t(t,x)+\h\dbH\big(t,x,\h\psi(t,x,\h V^{T'}_x(t,x),\h V^{T'}_{xx}(t,x)),\h V^{T'}_x(t,x),\h V^{T'}_{xx}(t,x)\big) =0,\q(t,x)\in[0,T]\times\dbR^n,\\
\ns\ds\h V^{T'}(T,x)\ges0=\h V^T(T,x),\qq x\in\dbR^n.\ea\right.\ee
For equations \rf{HJB4} and \rf{HJB'}, by the comparison theorem of (viscosity) solutions to parabolic PDEs, we obtain
\bel{V<V}0\les\h V^T(t,x)\les\h V^{T'}(t,x)\les K,\qq(t,x)\in[0,T]\times\dbR^n,\q T<T'<\i.\ee
Furthermore, for any fixed $(t,x)$ and admissible control $\h u(\cd)$,  it is obvious that
$$0\les\h J^{\,T'}(t,x;\h u(\cd))-\h J^{\,T}(t,x;\h u(\cd))=\dbE\int_T^{T'}\h g(s,\h X(s),\h u(s))ds\les\int_T^\i\h\f(s)ds.$$
Combing the definition of $\h V^T$, we get  \bel{V-V}0\les\h V^{T'}(t,x)-\h V^T(t,x)\les \int_T^\i\h\f(s)ds,\qq(t,x)\in[0,T]\times\dbR^n.\ee
From \rf{V<V} and \rf{V-V}, it follows that \rf{V to V} holds for some continuous function $\h V^\i(\cd\,,\cd)$,
$$\lim_{T\to\i}\h V^T(t,x)=\h V^\i(t,x),\qq(t,x)\in[0,\h T]\times\dbR^n,$$
uniformly, where $\h T\in\dbR_+$ is given.
Moreover, for any given $T\in\dbR_+$, by \rf{V-V}, we see that (sending $T'\to\i$)
$$0\les\h V^\i(t,x)-\h V^T(t,x)\les\int_T^\i\h\f(s)ds,\qq (t,x)\in[0,T]\times\dbR^n.$$
In particular, at $t=T$, we have
$$0\les\h V^\i(T,x)\les\int_T^\i\h\f(s)ds,\qq x\in\dbR^n,$$
which implies that
$\ds \lim_{T\to\i}\h V^\i(T,x)=0 $ uniformly in  $x\in\dbR^n$.

\ms

{\it Step 2: Convergence of the derivatives and regularity of the limits.} By the classical regularity theory for fully nonlinear parabolic equations (see \cite{Krylov 1987}), we get the following uniform estimate
$$\|\h V^T(\cd,\cd)\|_{C^{1+\frac\a 2,2+\a}([0,T]\times\dbR^n)}\les C,$$
for some $\a\in (0,1)$ and constant $C > 0$ independent of $T$. This implies that, for any compact set $\cK\subset\dbR_+\times\dbR^n$ and all sufficiently large $T$, the families $\{\h V_t^T(\cd,\cd)\}$, $\{\h V_x^T(\cd,\cd)\}$ and $\{\h V_{xx}^T(\cd,\cd)\}$ are equicontinuous on $\cK$.

\ms

Applying the Arzel\`{a}-Ascoli Theorem to these sequences and employing  a standard diagonal argument, we can  extract a common subsequence   $\{T_k\}_{k\ges 1}$ such that
$$\lim_{k\to\i}\sup_{\cK}\(|\h V_t^{T_k}(t,x)-\h V_t^\i(t,x)|+|\h V_x^{T_k} (t,x)-\h V_x^\i(t,x)|+|\h V_{xx}^{T_k}(t,x)-\h V_{xx}^\i(t,x)|\)=0.$$
Moreover, the limit function $\h V^\i(\cd,\cd)$ inherits the differentiability properties, i.e., $\h V^\i(\cd,\cd)\in C^{1,2}([0,\i) \times\dbR^n)$.

\ms

{\it Step 3: Verification of the limiting HJB equation.}
Consider the following finite horizon HJB equation,
\bel{HJB-Tk}\left\{\2n\ba{ll}
\ns\ds\h V^{T_k}_t(t,x)+\h\dbH \big(t,x,\widehat{\psi}(t,x,\widehat{V}^{T_k}_x(t,x),\widehat{V}^{T_k}_{xx}(t,x)), \widehat{V}^{T_k}_x(t,x),\widehat{V}^{T_k}_{xx}(t,x) \big) = 0, \q (t,x) \in[0,T_k]\times\dbR^n,\\
\ns\ds\h V^{T_k}(T_k,x)=0,\qq x\in\dbR^n.\ea\right.\ee
Then,  passing to the limit in the above equation, we  get
$$\ba{ll}
\ns\ds0=\lim_{k\to\i}\[\h V^{T_k}_t(t,x)+\h\dbH \big(t,x,\h\psi(t,x,\h V^{T_k}_x(t,x),\h V^{T_k}_{xx}(t,x)),\h V^{T_k}_x(t,x),\h V^{T_k}_{xx}(t,x)\big)\]\\
\ns\ds\q=\h V^\i_t(t,x)+\h\dbH\big(t,x,\h\psi(t,x,\h V^\i_x(t,x),\h V^\i _{xx}(t,x)),\h V^\i_x(t,x),\h V^\i_{xx}(t,x)\big),\ea$$
for all $(t,x)\in\dbR_+\times\dbR^n$. Combined with the terminal condition
$$\ds\lim_{T\to\i}\h V^\i(T,x)=0,\qq x\in\dbR^n,$$
we conclude that $\h V^\i(\cd,\cd)$ is a classical solution of the infinite horizon HJB equation \rf{HJB1}.

\ms

Finally, the uniqueness of the classical solution of \rf{HJB1} comes from the following Theorem \ref{Th-V-infty}. \endpf

\ms

Now, we can establish  the following result for  Problem ($\h{\rm C}$)$^\i$.

\bt{Th-V-infty} \sl Let {\rm($\h{\rm H}$1)--($\h{\rm H}$2)} and \rf{hpsi} hold. Then the value function $\h V^\i(\cd\,,\cd)$ of Problem {\rm($\h{\rm C}$)$^\i$} is the unique classical solution of the HJB equation \rf{HJB1}. Furthermore, for $(t,x)\in\dbR_+\times\dbR^n$, the open-loop optimal control of Problem $(\h{\rm C})^\i$ uniquely exists and is given by
\bel{hu^i}\h u^\i(s)=\h\psi\big(s,\h X^\i(s),\h V_x^\i(s,\h X^\i(s)),\h V_{xx}^\i(s,\h X^\i(s))\big),\qq s\in[t,\i),\ee
where $\h X^\i (\cd)$ is the unique strong solution to the following SDE
\bel{state-C-opt}\left\{\2n\ba{ll}
\ds d\h X^\i (s)=b\big(s,\h X^\i(s),\h\Psi(s,\h X^\i(s))\big)ds+\sum_{i=1}^d\si_i\big(s,\h X^\i(s),\h\Psi(s,\h X^\i(s))\big)dW_i(s),\qq s\in[t,\i),\\
\ns\ds\h X^\i(t)=x,\ea\right.\ee
with
$$\h\Psi(s,x)=\h\psi(s,x,\h V_x^\i(s,x),\h V_{xx}^\i(s,x)),\qq(s,x)\in\dbR_+\times\dbR^n.$$
Here $\h\Psi(\cd\,,\cd)$ is called the optimal feedback strategy of Problem {\rm($\h{\rm C}$)$^\i$}.

\et

\it Proof. \rm Let $\h V^\i(\cd\,,\cd)\in C^{1,2}(\dbR_+\times\dbR^n)$  be the solution to the HJB equation \rf{HJB1} and  $(\h u(\cd),\h X(\cd))$  be an admissible pair of \rf{hat X}. For any $T>t$, applying It\^o's formula to $s\mapsto\h V^\i(s,\h X(s))$ on $[t,T]$, we obtain
$$\ba{ll}
\ns\ds\dbE[\h V^\i(T,\h X(T))]-\h V^\i(t,x)\\
\ns\ds=\dbE\int_t^T\[\h V_s^\i(s,\h X(s))+\h\dbH\big(s,\h X (s),\h  u(s),\h V_x^\i(s,\h X(s)),\h V_{xx}^\i(s,\h X(s))\big)-\h g(s,\h X(s), \h u(s))\]ds\\
\ns\ds\ges-\dbE\int_t^T\h g(s,\h X(s),\h u(s))ds.\ea$$
Thus, for all $(t,x)\in\dbR_+\times\dbR^n$ and $\h u(\cd)\in\sU[t,\i]$,
$$\h V^\i(t,x)\les\dbE\[\int_t^T\h g(s,\h X(s),\h u(s))ds+\h V^\i(T,\h X(T))\].$$
Using the terminal condition $\ds\lim_{T\to\i}\h V^\i(T,x)=0$, $x\in\dbR^n$ and the boundedness estimate of $\h X(\cd)$, we conclude
$$\h V^\i(t,x)\les\dbE\int_t^\i\h g(s,\h X(s),\h u(s))ds=\h J^\i(t,x;\h u(\cd)),\q\forall\h u(\cd)\in\sU[t,\i].$$
By the arbitrariness of control $\h u(\cd)$, we have
$$\h V^\i(t,x)\les\inf_{\h u(\cd)\in\sU[t,\i]}\h
J^\i(t,x;\h u(\cd)),\q (t,x)\in[0,T]\times\dbR^n.$$
On the other hand, notice the definition of $\h\psi(\cd\,,\cd\,,\cd\,,\cd)$, we replace $\h u(\cd)$ in the above by
$$\h u^\i(\cd)=\h\psi\big(\cd\,,\h X^\i(\cd),\h V_x^\i(\cd\,,\h X^\i(\cd)),\h V_{xx}^\i(\cd\,,\h X^\i(\cd))\big),$$
where $\h X^\i(\cd)$ is the solution of SDE \rf{state-C-opt}. In this case, all the inequalities above turn  out to be equalities. Thus, we get
$$\h V^\i(t,x)=\dbE\int_t^\i\h g(s,\h X^\i(s),\h u^\i(s))ds= \h J^\i (t,x;\h u^\i(\cd)),\q (t,x)\in[0,T]\times\dbR^n.$$
This shows that $\h u^\i(\cd)$ is optimal and that $\h V^\i(\cd\,,\cd) $ is indeed the value function of Problem {\rm$(\h{\rm C})^\i$}. \endpf

\ms

We have seen that our approach is different from that in \cite{Baumeister-Leitao-Silva 2007}.

\section{Time-Inconsistent Problems in Infinite Horizon}

\rm

In this section, we will present a construction of equilibrium strategy for infinite horizon time-inconsistent optimal control problem, under condition \rf{l1}. We begin with the following hypotheses (compare with ($\h{\rm H}1$)--($\h{\rm H}2$)).

\ms

{\bf(H1)} The maps $b,\si_1,\cds,\si_d:\dbR_+\times\dbR^n\times U\to\dbR^n$ are continuous, regular enough, and for some $L>0$,
$$|b(t,0,u)|+\sum_{i=1}^d|\si_i(t,0,u)|\les L,\qq(t,u)\in\dbR_+\times U,$$
and
$$|b(t,x,u)-b(t,x',u)|+\sum_{i=1}^d|\si_i(t,x,u)-\si_i(t,x',u)|\les L|x-x'|,\q(t,u)\in\dbR_+\times U,~x,x'\in\dbR^n.$$
Moreover, the following non-degeneracy condition holds: For some $\e>0$,
\bel{sisi>d*}\si(t,x,u)\si(t,x,u)^\top\ges\e I_n,\q\forall(t,x,u)\in\dbR_+\times\dbR^n\times U.\ee

\ms

{\bf(H2)} The map $g:\dbR_+\times\dbR^n\times U\to\dbR_+$ is continuous, regular enough, and for some $\f(\cd)\in L^1(\dbR_+;\dbR)$, it holds that
\bel{g}0\les g(t,x,u)\les\f(t),\qq(t,x,u)\in\dbR_+\times\dbR^n\times U.\ee

\ms

Now, let $\t\in\dbR_+$ be given. From \rf{J=J+V}, we have
\bel{J^i}J^\i(\t,x;u(\cd))=\dbE\1n\int_\t^{\t+T_0}\3n\l(s-\t)g(s,X(s),u(s))ds\1n+\1n e^{\d\t}\h J^{\,\d,\i}(\t\1n+\1n T_0,X(\t\1n+\1n T_0);u(\cd+\t\1n+\1n T_0)),\ee
where (by setting $r=\t+T_0$)
$$\h J^{\,\d,\i}(r,x;u(\cd))=\dbE\int_r^\i e^{-\d s}g(s,X(s),u(s))ds,$$
with $(X(\cd),u(\cd))$ being the state-control pair of the state equation \rf{state1}. Thus, on $[\t+T_0,\i)$, it is natural to formulate the following optimal control problem.

\ms

{\bf Problem (C)$^{\t+T_0,\i}$.} For given $(t,x)\in[\t+T_0,\i)\times\dbR^n$, find $\bar u(\cd)\in\sU[t,\i]$ such that
$$\h J^{\,\d,\i}(t,x;\bar u(\cd))=\inf_{u(\cd)\in\sU[t,\i]}\h J^{\,\d,\i}(t,x;u(\cd))\equiv\h V^\d(t,x),$$
where the state-control pair satisfies \rf{state1} on $[t,\i)$.

\ms

By the result from Section 2, we have the following.

\bp{Pro-V-delta} \sl Let {\rm (H1)--(H2)} and \rf{l1} hold. Then the following admits the unique classical solution:
\bel{HJB2}\left\{\2n\ba{ll}
\ds V^{\d,\i}(t,x)+\dbH^\d\big(t,x,\psi^\d(t,x,V^{\d,\i}_x(t,x),     V^{\d,\i}_{xx}(t,x)),V^{\d,\i}_x(t,x),V^{\d,\i}_{xx}(t,x)\big)=0,\\
\ns\ds\qq\qq\qq\qq\qq\qq\qq\qq\qq\qq(t,x)\in[\t\1n+\1n T_0,\i)\times\dbR^n,\\
\ns\ds V^{\d,\i}(\i,x)=0,\qq x\in\dbR^n,\ea\right.\ee
where
\bel{psi1}\left\{\2n\ba{ll}
\ds\dbH^\d(t,x,u,p^ \top\2n,P)=p
b(t,x,u)+{1\over2}\sum_{i=1}^d\si_i(t,x,u)^ \top
P\si_i(t,x,u)+e^{-\d t}g(t,x,u),\\
\ns\ds\qq\qq\qq\qq\qq\qq\qq\qq(t,x,u,p,P)\in[\t+T_0,\i)\times\dbR^n\times U\times\dbR^n\times\dbS^n,\\
\ns\ds\psi^\d(t,x,p^\top\2n,P)\in\arg\min_{u\in U}\dbH^\d(t,x,u,p^ \top\2n,P),
\qq(t,x,p,P)\in[\t+T_0,\i)\times\dbR^n\times\dbR^n\times\dbS^n,\ea\right.\ee
with $\psi^\d(\cd\,,\cd\,,\cd\,,\cd)$ being regular enough. Further, Problem {\rm(C)$^{\t+T_0,\i}$} admits a unique open-loop optimal control $\bar u(\cd)$, and it has the following representation:
\bel{baru}\bar
u(s)=\psi^\d(s,\bar X(s),V^{\d,\i}_x(s,\bar X(s)),V^{\d,\i}_{xx}(s,\bar X(s))),\qq\q s\in[\t+T_0,\i),\ee
where $\bar X(\cd)$ is the corresponding open-loop optimal state process.
\ep

\ms

Next, on $[\t,\t+T_0]$, we have the state equation \rf{state1} with the following cost functional (see \rf{J=J+V}):
\bel{}J^\i(\t,x;u(\cd))=\dbE\[\int_\t^{\t+T_0}\l(s-\t)g(s,X(s),u(s))ds
+e^{\d\t}V^{\d,\i}(\t+T_0,X(\t+T_0))\].\ee
Then we have a time-inconsistent problem on a finite time horizon, called it Problem (N)$^{\t,\t+T_0}$, which can be solved following arguments of \cite{Yong 2012}. We now make it more precise. To this end, let us introduce the following definition.

\bde{steategy} Let $\t\in\dbR_+$ be given. A continuous map $\Psi:[\t,\i)\times\dbR^n\to U$ is called a {\it feedback strategy} for the system \rf{state1} on interval $[\t,\i)$, if for any $x\in\dbR^n$,
\bel{state-Psi}\left\{\2n\ba{ll}
\ds dX(s)=b\big(s,X(s),\Psi(s,X(s))\big)ds+\sum_{i=1}^d\si_i\big(s,X(s),
\Psi(s, X(s))\big)dW_i(s),\qq s\in[\t,\i),\\
\ns\ds X(\t)=x,\ea\right.\ee
admits a unique solution $X(\cd)\equiv X(\cd\,;\t,x,\Psi(\cd))$ on $[\t,\i)$.

\ede

\bde{equilibrium} Let $\t\in\dbR_+$ be given. A strategy $\Psi(\cd\,,\cd)$ on $[\t,\i)$ is called an {\it equilibrium strategy} of Problem (N)$^\i$ on $[\t,\i)$ if the following holds:

\ms

(i) The open-loop optimal control of Problem {\rm$(C)^{\t+T_0,\i}$} is given by
$$\bar u(s)=\Psi(s,\bar X(s)),\qq s\in[\t+T_0,\i).$$

(ii) $\Psi(\cd\,,\cd)$ is a time-consistent equilibrium strategy (in the sense of \citedef{Yong 2012}{Definition 4.1}) of Problem (N)$^{\t,\t+T_0}$.

\ede

\ms

When the above happens, we call $(\bar X(\cd),\bar u(\cd))$ a (time-consistent) equilibrium pair for the initial pair $(\t,x)$. Moreover, the equilibrium value at $(t,\bar X(t))$ is given by
$$\bar V^\i\big(t,\bar X(t;\t, x,\Psi(\cd))\big)=\dbE\[\int_t^{\t+T_0}\l(s-\t)g\big(s,\bar X(s),\Psi(s,\bar X(s))\big)ds+e^{\d\t}V^{\d,\i}\big(\t+T_0,\bar X(\t+T_0 )\big )\].$$

\ms

Now, let us summarize the above. Besides \rf{psi1}, we introduce the following.
\bel{psi2}\left\{\2n\ba{ll}
\ds\dbH^\l(\varrho ,t,x,u,p^ \top\2n,P)=pb(t,x,u)+{1\over2}\sum_{i=1}^d\si_i
(t,x,u)^ \top P\si_i(t,x,u)+\l(t-\varrho )g(t,x,u),\\
\ns\ds\psi^\l(\varrho,t,x,p^\top\2n,P)\in\arg\min_{u\in U}\dbH^\l(\varrho,t,x,u,p^\top\2n,P),\q  (\varrho,t,x,p,P)\in\D^*[0,T]\times\dbR^n\times\dbR^n\times\dbS^n.\ea\right.\ee

To emphasize the main idea, we still assume  $\psi^\l(\cd\,,\cd\,,\cd\,,\cd\,,\cd)$ to be regular enough.

\bt{Th-equilibrium} \sl Let {\rm (H1)--(H2)} and \rf{l1} hold.
%Then the following admit the unique classical solution:
%%
%\bel{HJB2}\left\{\2n\ba{ll}
%%
%\ds V^{\d,\i}(t,x)+\dbH^\d\big(t,x,\psi^\d(t,x,V^{\d,\i}_x(t,x),     V^{\d,\i}_{xx}(t,x)),V^{\d,\i}_x(t,x),V^{\d,\i}_{xx}(t,x)\big)=0,\\
%%
%\ns\ds\qq\qq\qq\qq\qq\qq\qq\qq\qq\qq(t,x)\in[\t\1n+\1n T_0,\i)\times\dbR^n,\\
%%
%\ns\ds V^{\d,\i}(\i,x)=0,\qq x\in\dbR^n.\ea\right.\ee
%%
Let the following admit a unique classical solution:
\bel{HJB3}\left\{\2n\ba{ll}
\ds\Th_t^\t(\varrho,t,x)+\dbH^\l\big(\varrho,t,x,\psi^\l(t,t,x,\Th_x^\t(t,t,x),
\Th_{xx}^\t(t,t,x)),\Th_x^\t(\varrho,t,x),\Th_{xx}^\t(\varrho,t,x)\big)=0,\\
\ns\ds\qq\qq\qq\qq\qq\qq\qq\qq\qq\qq\qq\t\les\varrho\les\t+T_0,~(t,x)\in[\varrho,\t+T_0]\times\dbR^n,\\
\ns\ds\Th^\t(\varrho,\t+T_0,x)=e^{\d\t}V^{\d,\i}(\t+T_0,x),\q x\in\dbR^n.\ea\right.\ee
Then the following gives an equilibrium value function of Problem {\rm(N)$^\i$}:
\bel{V}V^{\t,\i}(t,x)=\left\{\2n\ba{ll}
\ds\Th^\t(t,t,x),\qq\qq\qq (t,x)\in[\t,\t+T_0)\times\dbR^n,\\
\ns\ds e^{ \d \t}V^{\d,\i}(t,x),\qq\qq(t,x)\in[\t+T_0,\i)\times\dbR^n,\ea\right.\ee
and the following is the equilibrium strategy of Problem {\rm(N)$^\i$} at $\t\in\dbR_+$:
\bel{u}\ba{ll}
\ds\Psi(t,x)=\left\{\2n\ba{ll} \ns\ds \psi^\l\big(t,t,x,\Th_x^\t(t,t,x),\Th_{xx}^\t(t,t,x)\big),\qq (t,x)\in[\t,\t+T_0]\times\dbR^n,\\
\ns\ds\psi^\d(t,x,V^{\d,\i}_x(t,x),V^{\d,\i}_{xx}(t,x)),\qq\q\ (t,x)\in [\t+T_0,\i)\times\dbR^n,\ea\right.
\ea\ee
where $V^{\d,\i}(\cd,\cd)$ and $\psi^\d(\cd\,,\cd\,,\cd\,,\cd)$ are the ones in Proposition \ref{Pro-V-delta}.
\et

We note that HJB equation \rf{HJB2} is a Cauchy problem for a non-degenerate PDE, which admits a unique solution $V^{\d,\i}(\cd\,,\cd)$. For system \rf{HJB3}, when $\si(t,x,u)=\si(t,x)$ (independent of $u$), we have a unique classical (global) solution $\Th^\t(\cd\,,\cd\,,\cd)$ (see \cite{Yong 2012}). In general, under suitable conditions, we have a unique classical local solution (see \cite {Lei-Pun 2023}). It remains open for the global solvability of \rf{HJB3}.

\ms

\ms

We emphasize that if $\t\in\dbR_+$ is current time, by finding an equilibrium strategy $\Psi(\cd\,,\cd)$, we mean that the applied open-loop control $\Psi(s,\bar X(s))$ on $[\t,\i)$ will remain on the time interval, keeping the local optimality, where $\bar X(\cd)=X(\cd\,;\t,x,\Psi(\cd))$ is the corresponding state process.

\ms

Now, we look at some possible extension of the above. Let $t\in\dbR_+$ be given. We consider a more general cost functional
\bel{ext}\ba{ll}
\ns\ds J^\i(t,x;u(\cd))=\dbE\int_t^\i g(t,s,X(s),u(s))ds\\
\ns\ds\qq\qq\qq=\dbE\[\int_t^{t+T_0}g(t,s,X(s),u(s))ds+\int_{t+T_0}^\i e^{-\d(s-t)}g_0(s,X(s),u(s))ds\]\\
\ns\ds\qq\qq\qq=\dbE\[\int_t^{t+T_0}g(t,s,X(s),u(s))ds+e^{\d t}\int_{t+T_0}^\i e^{-\d s}g_0(s,X(s),u(s))ds\],\ea\ee
which is true if instead of \rf{l1}, we assume the following:
\bel{l2}g(t,s,x,u)=e^{-\d(s-t)}g_0(s,x,u),\q(t,s,x,u)\in\D^*[0,\i]\times\dbR^n\times U,~s-t\ges T_0.\ee
Note that on the right-hand side of \rf{ext}, the second integral is on $[t+T_0,\i)$, which implies
$$s-t\ges T_0.$$
Thus, \rf{ext} holds if \rf{l2} is assumed. Hence, in this case, for any $t\in\dbR_+$ given, divided the interval $[t,\i)$ into two parts: $[t,t+T_0]$ and $[t+T_0,\i)$. On the second one, the corresponding optimal control problem is a classical, and on the first one, it is a time-inconsistent optimal control problem (in a finite interval). Then, we can find an equilibrium strategy for the problem, using our main idea with some small modifications.

\section{Recursive Cost Functional Case.}

In this section, we consider the case of the so-called recursive cost functional. As we know that (see \cite{Wei-Yong-Yu 2017}) the recursive functional in a finite time horizon $[t,T]$ takes the following form
\bel{J}J^T_r(t,x;u(\cd))=Y(t;t,x),\ee
for the state-control pair $(X(\cd),u(\cd))$ of, say, the state equation \rf{state1}, where $(Y(\cd\,;t,x),Z(\cd\,;t,x))$ is the adapted solution of the following BSDE:
\bel{BSDE1}\left\{\2n\ba{ll}
\ns\ds dY(s)=-g(t,s,X(s),u(s),Y(s),Z(s))ds+Z(s)dW(s),\qq s\in[t,T],\\
\ns\ds Y(T)=h(t,X(T)),\ea\right.\ee
with $t$ being a parameter and with some $g:\D^*[0,T]\times\dbR^n\times U\times\dbR\times\dbR^{1\times d}\to\dbR$, $h:[0,T]\times\dbR^n\to\dbR$ (\cite{Lazrak-Quenez 2003,Lazrak 2004}). By looking at the integral form of \rf{BSDE1}, we see that \rf{J} can be written as
\bel{J1}J^T_r(t,x;u(\cd))=\dbE\[\int_t^Tg(t,s,X(s),u(s),Y(s),Z(s))ds
+h(t,X(T))\].\ee
Now, we introduce the following cost functional
\bel{J2}J^\i_r(t,x;u(\cd))=\dbE\int_t^\i g(t,s,X(s),u(s),Y(s),Z(s))ds,\ee
assuming the following
$$0\les g(t,s,x,u,y,z)\les\f(s),\q(t,s,x,u,y,z)\in\D^*[0,\i]\times\dbR^n\times U\times\dbR\times\dbR^{1\times d},$$
with $\f(\cd)\in L^1(\dbR_+;\dbR)$. Thus, the cost functional is well-defined. Then we formulate the following problem.

\ms

{\bf Problem (P)$_r^\i$.} For given $t\in\dbR_+$, find $\bar u(\cd)\in\sU[t,\i]$, such that
$$J^\i_r(t,x;\bar u(\cd))=\inf_{u(\cd)\in\sU[t,\i]}J^\i_r(t,x;u(\cd)).$$

\ms

In general, the above problem is time-inconsistent (\cite{Wei-Yong-Yu 2017}). Suggested by \rf{l2}, we introduce the following.

\ms

{\bf(H3)} The maps $g:\D^*[0,\i]\times\dbR^n\times U\times\dbR\times\dbR^{1\times d}\to\dbR_+$ and $g_0:\dbR_+\times\dbR^n\times U\to\dbR_+$ are continuous, bounded, and
for some $\d>0$, it holds
\bel{g0*}g(t,s,x,u,y,z)=-\d y+g_0(s,x,u),\qq(t,s,x,u,y,z)\in\D^*\times\dbR^n\times U\times\dbR\times\dbR^{1\times d},~s-t\ges T_0.\ee

\ms

Then under (H3), we have
\bel{J3}J^\i_r(t,x;u(\cd))=\dbE\[\int_t^{t+T_0}g(t,s,X(s),u(s),Y(s),Z(s))ds
+\int_{t+T_0}^\i\(-\d Y(s)+g_0(s,X(s),u(s))\)ds\].\ee
Now, we introduce the following BSDE:
\bel{BSDE}\left\{\2n\ba{ll}
\ns\ds d\wt Y(s)=\(\d\wt Y(s)-g_0(s,X(s),u(s))\)ds+\wt Z(s)dW(s),\q s\in[\t,\i),\\
\ns\ds\wt Y(\i)=0.\ea\right.\ee
Then, on one hand,
$$\wt Y(t+T_0)=\dbE_{t+T_0}\int_{t+T_0}^\i\(-\d\wt Y(s)+g_0(s,X(s),u(s))\)ds.$$
On the other hand, by variation of constants formula, we have
\bel{wtY}\wt Y(t+T_0)=\dbE_{t+T_0}\int_{t+T_0}^\i e^{-\d(s-\t)}g_0(s,X(s),u(s))ds,\qq t\in[\t,\i).\ee
Thus, by uniqueness, one has
\bel{J4}J^\i_r(t,x;u(\cd))=\dbE\[\int_t^{t+T_0}g(t,s,X(s),u(s),Y(s),Z(s))ds
+\int_{t+T_0}^\i e^{\d(s-t-T_0)}g_0(s,X(s),u(s))ds\].\ee
We know from the introduction that if the cost functional is of the second integral, the corresponding optimal control problem is time-consistent. Hence, we divide the integral $[t,\i)$ into two parts: $[t,t+T_0]$ and $[t+T_0,\i)$. In the second one, we solve a classical optimal control problem as in Section 3, to get its value function $V_r(\cd\,,\cd)$. Then we have a naive optimal control problem in a finite interval with the following cost functional:
\bel{J_r}J_r^\i(t,x;u(\cd))=\dbE\[\int_t^{\t+T_0}g(t,s,X(s),u(s),Y(s),Z(s))ds+
e^{\d(t+T_0)}V_r(t+T_0,X(t+T_0))\].\ee
Equivalently,
\bel{J_r}J^\i_r(t,x;u(\cd))=Y(t;t,x),\ee
where $(Y\cd),Z(\cd))$ is the adapted solution of the following BSDE:
\bel{BSDE2}\left\{\2n\ba{ll}
\ns\ds dY(s)=-g(t,s,X(s),u(s),Y(s),Z(s))ds+Z(s)dW(s),\qq s\in[t,t+T_0],\\
\ns\ds Y(t+T_0)=e^{\d(t+T_0)}V_r(t+T_0,X(t+T_0)),\ea\right.\ee
We now pose the following problem.

\ms
{\bf Problem (N)$^\i_r$.} \rm For any initial pair $(t,x)\in[0,\i)\times\dbR^n$, find a $\bar u(\cd\,;t,x)\in\sU[t,\i]$ such that
\bel{V2*}J^\i_r(t,x;\bar
u(\cd\,;t,x))=\inf_{u(\cd)\in\sU[t,\i]}J^\i_r(t,x;u(\cd))\equiv V^\i_r(t,x).\ee

\ms

Similar to before, this is a naive problem. Such a naive problem can be solved by combining the main idea of the previous section and \cite{Wei-Yong-Yu 2017}. We omit the details here.

\ms

To conclude this section, we look at possible extensions.

\ms

From the above, we have seen that for given $t\in\dbR_+$, we can naturally divide the interval $[t,\i)$ into two parts: $[t,t+T_0]$ and $[t+T_0,\i)$ for some $T_0>0$. Now if it turns out that on $[t+T_0,\i)$, the corresponding optimal control problem is time-consistent, then we can solve such a problem first to get the value function, denoted by $V_r(\cd\,,\cd)$. Hence, we may consider the state equation \rf{state1} with cost functional \rf{J_r} with $(Y(\cd),Z(\cd))$ being the adapted solution to the BSDE:
\bel{BSDE3}\left\{\2n\ba{ll}
\ns\ds dY(s)=-g(t,s,X(s),u(s),Y(s),Z(s))ds+Z(s)dW(s),\qq s\in[t,\i),\\
\ns\ds Y(\i)=0.\ea\right.\ee
We now introduce the following hypothesis.

\ms

{\bf(H4)} Maps $g:\D^*[0,\i)\times\dbR^n\times U\times\dbR\times\dbR^{1\times d}\to\dbR_+$ and $g_0:\dbR_+\times\dbR^n\times U\times\dbR\times\dbR^{1\times d}\to\dbR_+$ are continuous. There exists a $T_0>0$ such that
\bel{g*}g(t,s,x,u,y,z)=g_0(s,x,u,y,z),\qq(t,s,x,u,y,z)\in\D^*[0,\i)\times\dbR^n\times U\times\dbR\times\dbR^{1\times d},~s-t\ges T_0.\ee
Further, for some $\f(\cd)\in L^1(\dbR_+;\dbR)$, it holds
\bel{les f}|g_0(s,x,u,y,z)|\les\f(s),\q(s,x,u,y,z)\in\dbR_+\times\dbR^n\times U\times\dbR\times\dbR^{1\times d}.\ee

\ms

Under (H4), for any $t\in\dbR_+$ given, we have
$$\ba{ll}
\ns\ds J^\i(t,x;u(\cd))=\dbE\int_t^\i g(t,s,X(s),u(s),Y(s),Z(s))s\\
\ns\ds\qq=\dbE\[\int_t^{t+T_0}g(t,s,X(s),u(s),Y(s),Z(s))ds+\int_{t+T_0}^\i
g_0(s,X(s),u(s),Y(s),Z(s))ds\].\ea$$
Thus, under \rf{les f}, it is easy for us to formulate a time-consistent optimal control problem with a recursive cost functional on $[t+T_0,\i)$. Then using the arguments of \cite{Wei-Yong-Yu 2017}, together with the idea of the previous section, we can also get the equilibrium strategy for the related naive problem. We leave the details to the interested readers.

\section{Conclusions.}

\ms

In the current paper, we have discussed the time-inconsistent optimal control problems in an infinite time horizon. The crucial assumption is \rf{l1} (or its extended version \rf{l2}). Such a condition should be reasonable for which most normal people will have such a characterization. Under such a condition, we are able to construct an equilibrium strategy. Further, similar idea has been used to the case of recursive cost functional problem. For such problems, instead of \rf{l1}, we have assumed \rf{g0*} (or extended version \rf{g*}) and obtained a similar result.

\ms

It is left widely open that what happens if the non-degenerate condition \rf{sisi>d*} (or \rf{sisi>d}) is not assumed. For that, we need to have the vanishing of viscosity to obtain the existence some sort of the solution to the equilibrium HJB equation. We are not able to get that for the time being.

\ms

Another open question is the global well-posedness is the equilibrium HJB equation \rf{HJB3}. The local well-posedness of them were obtained in \cite{Lei-Pun 2023}. We suspect that if one uses some kind of monotonicity conditions, one might be able to get the expected results.

\end{document}